\documentclass{daj}
\usepackage{amsmath, amsfonts, amssymb}

\global\long\def\tr{\mathrm{{tr}}}%

\global\long\def\a{\alpha}%

\global\long\def\b{\beta}%

\global\long\def\E{\mathrm{\mathbb{E}}}%

\global\long\def\P{\mathrm{\mathbb{P}}}%

\dajAUTHORdetails{%
title = {Quasirandom groups enjoy interleaved mixing}, %% please capitalize all significant words
author = {Harm Derksen and Emanuele Viola},
plaintextauthor = {Harm Derksen and Emanuele Viola},
plaintexttitle = {Quasirandom groups enjoy interleaved mixing}, %% title without math or LaTeX
keywords = {representation theory, mixing, group, interleaved},
} %%% END \dajAUTHORdetails

\dajEDITORdetails{%
year={2023},
number={14},
received={23 June 2022}, % received date: example: 7 January 2017
published={4 September 2023}, % published date
doi={10.19086/da.84738}, % XXX = number of paper, e.g. da006 for paper#6
} %%% END \dajEDITORdetails

\begin{document}

\begin{frontmatter}[classification=text]
%% EDITOR: this will force the keywords to appear right after the Abstract.
%% If the abstract is too long and would force the keywords off the
%% front page, please comment out % [classification=text] above
%% This way the keywords will be floated on the bottom of the first page
%% even though the Abstract spills over to the next page.

%%% AUTHOR: Title goes here. This line is optional. You must use it
%% if title has footnote attached or requires nontrivial typesetting,
%% e.g., inclusion of linebreaks to force nice layout.
\title{Quasirandom groups enjoy interleaved mixing} %% please capitalize all significant words

%%% AUTHOR:
%%% List all authors. If you wish, place grant acknowledgements in \thanks.
%%% In brackets include a short tag for each author.
\author[HD]{Harm Derksen\thanks{Partially supported by NSF grant DMS 2147769.}}
\author[EV]{Emanuele
Viola\thanks{Supported by NSF grant CCF-2114116.}}

%%% AUTHOR: Abstract goes here
\begin{abstract}
Let $G$ be a group such that any non-trivial representation has dimension
at least $d$. Let $X=(X_{1},X_{2},\ldots,X_{t})$ and $Y=(Y_{1},Y_{2},\ldots,Y_{t})$
be distributions over $G^{t}$. Suppose that $X$ is independent from
$Y$. We show that for any $g\in G$ we have
\[
\left|\P[X_{1}Y_{1}X_{2}Y_{2}\cdots X_{t}Y_{t}=g]-1/|G|\right|\le\frac{|G|^{2t-1}}{d^{t-1}}\sqrt{\E_{h\in G^{t}}X(h)^{2}}\sqrt{\E_{h\in G^{t}}Y(h)^{2}}.
\]
Our results generalize, improve, and simplify previous works.

\end{abstract}
\end{frontmatter}

%%% AUTHOR: body of paper starts here
\bigskip

\emph{Quasirandom groups}, introduced by Gowers \cite{Gowers08},
are groups whose non-trivial representations have large dimension.
Multiplication in such groups is known to behave like a random function
in several respects. The prime example of this is that if $X$ and
$Y$ are independent, high-entropy distributions over a quasirandom
group then $XY$ (i.e., sample from each and output the product) becomes
closer to uniform in $L_{2}$ norm. For a discussion of this result
and its many proofs we refer to Section 13 of \cite{MR3584096}. Other
random-like behaviors are known with respect to, for example, \emph{progressions
}\cite{DBLP:conf/innovations/BhangaleHR22} and \emph{corners }\cite{Austin2016}
(cf.~\cite{viola-SIGACT19}).

In this work we are interested in a question posed by Miles and Viola
\cite{MilesV-leak}. Let $X=(X_{1},X_{2})$ and $Y=(Y_{1},Y_{2})$
be high-entropy distributions over $G^{2}$ such that $X$ is independent
from $Y$ (but $X_{1}$ needs not be independent from $X_{2}$ and
$Y_{1}$ needs not be independent from $Y_{2}$). They asked if the
\emph{interleaved product} $X_{1}Y_{1}X_{2}Y_{2}$ ``mixes,'' i.e.,
if it is close to uniform, for suitable groups $G$. Their question
was motivated by an application to cryptography (which follows from
a positive answer to a more general question they asked).

Gowers and Viola \nocite{GowersV-cc-int,GowersV-cc-int-2,GowersV-cc-int-journal}
give a positive answer to this question for non-abelian simple groups,
which are known to be quasirandom. For the special case of $G=SL(2,q)$
they prove a strong error bound. A simpler exposition of the latter
proof appears in \cite{viola-SIGACT19}. A follow-up paper by Shalev
\cite{Shalev16} gives stronger error bounds for non-abelian simple
groups.

These proofs are somewhat complicated and use substantial machinery,
and they only apply to simple groups. Here we give a very short and
elementary proof that applies to any quasirandom group, as stated
in the abstract.

To illustrate the bound in the abstract, suppose that $X$ is uniform
over a set of density $\a$ and $Y$ is uniform over a set of density
$\b$. Then the right-hand side is $|G|^{2t-1}\cdot d^{-t+1}\cdot(\a\b)^{-1/2}/|G|^{2t}=|G|^{-1}\cdot d^{-t+1}\cdot(\a\b)^{-1/2}.$
Our results also slightly improve the parameters in the cases where
interleaved mixing could be established. For example for $t>2$ the
bounds in \cite{GowersV-cc-int-journal} and \cite{Shalev16} have
$(\a\b)^{-1}$ instead of $(\a\b)^{-1/2}$.

The paper \cite{GowersV-cc-int-journal} also shows that from interleaved
mixing there follow a number of other results (including the solution
to the more general question in \cite{MilesV-leak}, thus enabling
the motivating application). Hence our results yield these applications
for any quasirandom group. Since this is an immediate composition
of proofs in \cite{GowersV-cc-int-journal} and this paper, we refer
the reader to \cite{GowersV-cc-int-journal} for precise statements.

\paragraph{Proof of the statement in the abstract.}

We follow standard notation for non-abelian Fourier analysis, see
for example Section 13 of \cite{MR3584096} or \cite{GowersV-group-mix}.
It suffices to prove the theorem for $g=1_{G}$. Let $Z$ be a distribution
over $G$. By Fourier inversion, and using that $\rho(1)=I$ and $\widehat{Z}(1)=1/|G|$
we have
\begin{equation}
|\P[Z=1]-1/|G||=|\sum_{\rho}d_{\rho}\tr(\widehat{Z}(\rho)\rho(1)^{T})-1/|G||=|\sum_{\rho\ne1}d_{\rho}\tr(\widehat{Z}(\rho))|\le\sum_{\rho\ne1}d_{\rho}|\tr(\widehat{Z}(\rho))|,\label{eq:error}
\end{equation}
where $\rho$ ranges over irreducible representations.

The main claim is that if $Z$ is the interleaved product $X_{1}Y_{1}X_{2}Y_{2}\cdots X_{t}Y_{t}$
where $X$ and $Y$ are as in the abstract then for any $\rho$
\begin{equation}
|\tr(\hat{Z}(\rho))|\le|G|^{2t-1}|\hat{X}(\rho^{\otimes t})|_{2}|\hat{Y}(\rho^{\otimes t})|_{2}.\label{eq:main-claim}
\end{equation}

Assuming the claim the proof is completed as follows. Plugging Inequality
\eqref{eq:main-claim} into \eqref{eq:error} and multiplying by $(d_{\rho}/d)^{t-1}$
which is $\ge1$ for $\rho\ne1$, the error is at most
\[
\frac{|G|^{2t-1}}{d^{t-1}}\sum_{\rho\ne1}\left(d_{\rho}^{t/2}\left|\hat{X}(\rho^{\otimes t})\right|_{2}\right)\left(d_{\rho}^{t/2}\left|\hat{Y}(\rho^{\otimes t})\right|_{2}\right).
\]

By Cauchy-Schwarz this is at most
\[
\frac{|G|^{2t-1}}{d^{t-1}}\sqrt{\sum_{\rho\ne1}d_{\rho}^{t}\left|\hat{X}(\rho^{\otimes t})\right|_{2}^{2}}\sqrt{\sum_{\rho\ne1}d_{\rho}^{t}\left|\hat{Y}(\rho^{\otimes t})\right|_{2}^{2}}.
\]

Note that $d_{\rho}^{t}$ is the dimension of $\rho^{\otimes t}$.

The representations $\rho^{\otimes t}$ are irreducible representations
of $G^{t}$, so each sum can be bounded above by summing over all
irreducible representations of $G^{t}$. Hence by Parseval the sum
with $X$ is at most $\E_{h\in G^{t}}X^{2}(h)$ and the same for $Y$,
proving the theorem.

Next we verify Inequality \eqref{eq:main-claim}. By definition we
have
\[
\hat{Z}(\rho)=\E_{g}Z(g)\overline{\rho(g)}=\E_{g}\sum_{g_{1},g_{2},\ldots,g_{2t}:\prod g_{i}=g}X(g_{1},g_{3},\ldots,g_{2t-1})Y(g_{2},g_{4},\ldots,g_{2t})\overline{\rho(g)}.
\]
This summation is the same as summing over all $g_{i}$ and setting
$g$ to be the product. Further, because $\rho$ is a representation
one has $\rho(\prod_{i}g_{i})=\prod_{i}\rho(g_{i})$. Hence we get
\[
\hat{Z}(\rho)=\frac{1}{|G|}\sum_{g_{1},g_{2},\ldots,g_{2t}}X(g_{1},g_{3},\ldots,g_{2t-1})Y(g_{2},g_{4},\ldots,g_{2t})\overline{\prod_{i\le2t}\rho(g_{i})}.
\]

And now the critical equation:
\begin{align*}
\tr\hat{Z}(\rho) & =\sum_{i}\frac{1}{|G|}\sum_{g_{1},g_{2},\ldots,g_{2t}}X(g_{1},g_{3},\ldots,g_{2t-1})Y(g_{2},g_{4},\ldots,g_{2t})\sum_{i_{2},i_{3},\ldots,i_{2t}}\bar{\rho}(g_{1})_{i,i_{2}}\bar{\rho}(g_{2})_{i_{2},i_{3}}\cdots\bar{\rho}(g_{2t})_{i_{2t},i}\\
& =\frac{1}{|G|}\sum_{i,i_{2},i_{3},\ldots,i_{2t}}\left(\sum_{g_{1},g_{3},\ldots,g_{2t-1}}X(g_{1},g_{3},\ldots,g_{2t-1})\bar{\rho}(g_{1})_{i,i_{2}}\cdot\bar{\rho}(g_{3})_{i_{3},i_{4}}\cdots\bar{\rho}(g_{2t-1})_{i_{2t-1},i_{2t}}\right)\\
& \ \cdot\left(\sum_{g_{2},g_{4},\ldots,g_{2t}}Y(g_{2},g_{4},\ldots,g_{2t})\bar{\rho}(g_{2})_{i_{2},i_{3}}\cdot\bar{\rho}(g_{4})_{i_{4},i_{5}}\cdots\bar{\rho}(g_{2t})_{i_{2t},i}\right)\\
& =|G|^{2t-1}\sum_{i,i_{2},i_{3},\ldots,i_{2t}}\left(\hat{X}(\rho^{\otimes t})_{i,i_{2},i_{3},\ldots,i_{2t}}\right)\left(\hat{Y}(\rho^{\otimes t})_{i_{2},i_{3},\ldots,i_{2t},i}\right).
\end{align*}
Inequality \eqref{eq:main-claim} now follows by applying the Cauchy-Schwarz
inequality.

%%% AUTHOR:
%%% Bibliography goes here. Note that the arXiv cannot process bibtex
%%% or biber bibliographies. Example of acceptable bibliograpy format:
\bibliographystyle{amsplain}

\begin{thebibliography}{99}
\bibitem[Aus16]{Austin2016}
Tim Austin.
\newblock Ajtai-{S}zemer\'edi theorems over quasirandom groups.
\newblock In {\em Recent trends in combinatorics}, volume 159 of {\em IMA Vol.
Math. Appl.}, pages 453--484. Springer, [Cham], 2016.

\bibitem[BHR22]{DBLP:conf/innovations/BhangaleHR22}
Amey Bhangale, Prahladh Harsha, and Sourya Roy.
\newblock Mixing of 3-term progressions in quasirandom groups.
\newblock In Mark Braverman, editor, {\em ACM Innovations in Theoretical
Computer Science conf.~(ITCS)}, volume 215 of {\em LIPIcs}, pages 20:1--20:9.
Schloss Dagstuhl - Leibniz-Zentrum f{\"{u}}r Informatik, 2022.

\bibitem[Gow08]{Gowers08}
W.~T. Gowers.
\newblock Quasirandom groups.
\newblock {\em Combinatorics, Probability {\&} Computing}, 17(3):363--387,
2008.

\bibitem[Gow17]{MR3584096}
W.~T. Gowers.
\newblock Generalizations of {F}ourier analysis, and how to apply them.
\newblock {\em Bull. Amer. Math. Soc. (N.S.)}, 54(1):1--44, 2017.

\bibitem[GV]{GowersV-cc-int-2}
W.~T. Gowers and Emanuele Viola.
\newblock The multiparty communication complexity of interleaved group
products.
\newblock {\em SIAM J.~on Computing}.

\bibitem[GV15]{GowersV-cc-int}
W.~T. Gowers and Emanuele Viola.
\newblock The communication complexity of interleaved group products.
\newblock In {\em ACM Symp.~on the Theory of Computing (STOC)}, 2015.

\bibitem[GV19]{GowersV-cc-int-journal}
W.~T. Gowers and Emanuele Viola.
\newblock Interleaved group products.
\newblock {\em SIAM J.~on Computing}, 48(3):554--580, 2019.
\newblock Special issue of FOCS 2016.

\bibitem[GV22]{GowersV-group-mix}
W.~T. Gowers and Emanuele Viola.
\newblock Mixing in non-quasirandom groups.
\newblock In {\em ACM Innovations in Theoretical Computer Science
conf.~(ITCS)}, 2022.

\bibitem[MV13]{MilesV-leak}
Eric Miles and Emanuele Viola.
\newblock Shielding circuits with groups.
\newblock In {\em ACM Symp.~on the Theory of Computing (STOC)}, 2013.

\bibitem[Sha16]{Shalev16}
Aner Shalev.
\newblock Mixing, communication complexity and conjectures of {G}owers and
{V}iola.
\newblock {\em Combinatorics, Probability and Computing}, pages 1--13, 6 2016.
\newblock arXiv:1601.00795.

\bibitem[Vio19]{viola-SIGACT19}
Emanuele Viola.
\newblock Non-abelian combinatorics and communication complexity.
\newblock {\em SIGACT News, Complexity Theory Column}, 50(3), 2019.

\end{thebibliography}

%% AUTHOR: You can generate such a bibliography from a .bib file by
%% running pdflatex/bibtex/pdflatex/pdflatex and then pasting the .bbl file
%% between \begin{thebibliography} and \end{bibliography}

%%% AUTHOR: Include a short description of each author following the
%%% structure below. Use the same short tags used previously.
%%% Use \imageat{} and \imagedot{} instead of "@" and "." in
%%% email addresses-this replaces the symbols with graphics to avoid
%%% e-mail address harvesting from the .pdf file
\begin{dajauthors}
\begin{authorinfo}[HD]
Harm Derksen\\
Northeastern University\\
Boston, MA\\
\url{https://hderksen.sites.northeastern.edu/}
\end{authorinfo}
\begin{authorinfo}[EV]
Emanuele Viola\\
Northeastern University\\
\url{https://www.ccs.neu.edu/home/viola/}
\end{authorinfo}
\end{dajauthors}

\end{document}